 \documentclass[11pt]{article}
 \usepackage{mathrsfs}
 \usepackage{amsfonts}

 \usepackage{mathrsfs,amsfonts}
 \usepackage{harvard}

 \newtheorem{theorem}{Theorem}[section]
 \newtheorem{lemma}[theorem]{Lemma}
 \newtheorem{corol}[theorem]{Corollary}
 \newtheorem{prop}[theorem]{Proposition}

 \def\blemma{\begin{lemma}}\def\elemma{\end{lemma}}
 \def\bproposition{\begin{prop}}\def\eproposition{\end{prop}}
 \def\btheorem{\begin{theorem}}\def\etheorem{\end{theorem}}
 \def\bcorollary{\begin{corol}}\def\ecorollary{\end{corol}}

 \def\beqlb{\begin{eqnarray}}\def\eeqlb{\end{eqnarray}}
 \def\beqnn{\begin{eqnarray*}}\def\eeqnn{\end{eqnarray*}}

 \def\ar{\!\!\!&}

 \def\d{\mbox{\rm d}}

 \def\proof{\noindent{\it Proof.~~}}\def\qed{\hfill$\Box$\medskip}

 \def\<{\langle}\def\>{\rangle}

 \def\mbb{\mathbb}\def\mbf{\mathbf}

\begin{document}

\noindent{(Version: 2012/4/11)}

\bigskip\bigskip

\centerline{\huge\bf Limit theorems for continuous}

\medskip

\centerline{\huge\bf time branching flows\footnote{Supported by NSFC and
985 Project.}}

\bigskip

\centerline{Hui He and Rugang Ma\footnote{Corresponding author. E-mail:
marugang@mail.bnu.edu.cn.}}

\bigskip

\centerline{Beijing Normal University}

\bigskip

{\narrower

\noindent\textit{Abstract:} We construct a flow of continuous time and
discrete state branching processes. Some scaling limit theorems for the
flow are proved, which lead to the path-valued branching processes and
nonlocal branching superprocesses over the positive half line studied in
Li (2012).

\smallskip

\noindent\textit{Mathematics Subject Classification (2010):} Primary
60J68, 60J80; secondary 60G57

\smallskip

\noindent\textit{Key words and phrases:} Stochastic flow, branching
process, continuous-time, discrete-state, superprocess, nonlocal
branching.

\par}

\bigskip


\section{Introduction}

A genealogical tree is naturally associated with a Galton-Watson branching
process. A continuous-state branching process (CB-process) can be obtained
as the small particle limit of rescaled Galton-Watson processes; see,
e.g., Lamperti (1967). The genealogical structures of binary branching
CB-processes were investigated by introducing \textit{continuum random
trees} in the pioneer work of Aldous (1991, 1993). Continuum random trees
corresponding to general branching mechanisms were constructed in Le~Gall
and Le~Jan (1998a, 1998b) and were studied further in Duquesne and Le~Gall
(2002). By pruning a Galton-Watson tree, Aldous and Pitman (1998) and Abraham at al. (2011)
constructed a tree-valued Markov process. Tree-valued processes associated
with general CB-processes were studied in Abraham and Delmas (2010) by
pruning arguments.

Motivated by the study of genealogy trees for critical branching processes
conditioned on non-extinction, Bakhtin (2011) studied a flow of binary
branching continuous-state branching processes with immigration
(CBI-processes) driven by a time-space Gaussian white noise. He also
pointed out the connection of the model with a superprocess conditioned on
non-extinction. In Li (2012), a class of path-valued branching processes
were constructed and studied using the techniques of stochastic equations
and superprocesses. The work is closely related to those of Bertoin and Le
Gall (2006) and Dawson and Li (2012). In a special case, the path-valued
branching processes in Li (2012) can be coded by the tree-valued processes
of Abraham and Delmas (2010). In He and Ma (2012), two flows of discrete
time and state Galton-Watson branching processes were introduced. There it
was showed that suitable rescaled sequences of those flows converge to
special forms of the flows of Dawson and Li (2012) and Li (2012),
respectively. The limit theorems in He and Ma (2012) were given in the
setting of the corresponding superprocesses. From those limit theorems the
convergence of the finite-dimensional distributions of corresponding the
path-valued processes was derived. The results give a better understanding
of the connection between discrete and continuum tree-valued branching
processes.

In this paper, we introduce a kind of flows of continuous time and
discrete state branching processes. We shall prove the scaling limit
theorems for those flows of the type of He and Ma (2012). In Section 2 a
short review is given to the path-valued branching processes and nonlocal
branching superprocesses studied in Li (2012). In Section 3 we construct a
continuous time and discrete state branching processes as the strong
solution of a stochastic integral equation. In Section 4 the construction
is extended to branching flows by considering stochastic equation systems.
In Section 5 we prove that suitable rescaled sequences of those flows
converge to the nonlocal branching superprocess. From the limit theorem we
also derive the convergence of the finite-dimensional distributions of
corresponding the path-valued processes.

Let $\mbb N$ = $\{0,1,2,\cdots\}$ and $\mbb N_+$= $\{1,2,\cdots\}$. Let $M[0,1]$ be the set of finite Borel measures on $[0,1]$
endowed with the topology of weak convergence. We identify $M[0,1]$ with
the set $F[0,1]$ of positive right continuous increasing functions on
$[0,1]$. Let $B[0,1]$ be the Banach space of bounded Borel functions on
$[0,1]$ endowed with the supremum norm $\|\cdot\|$. Let $C[0,1]$ denote
its subspace of continuous functions. We use $B[0,1]^+$ and $C[0,1]^+$ to
denote the subclasses of positive elements and $C[0,1]^{++}$ to denote the
subset of $C[0,1]^+$ of functions bounded away from zero. For $\mu\in
M[0,1]$ and $f\in B[0,1]$ write $\<\mu, f\>$ = $\int f d\mu$ if the
integral exists. Let $D([0,\infty),M[0,1])$ denote the space of c\`{a}dl\`{a}g paths from
$[0,\infty)$ to $M[0,1]$ endowed with the Skorokhod topology.
Throughout the paper, we only consider \textit{continuous
time} processes, so we shall often omit this phrase in the sequel.


\section{Preliminaries}

\setcounter{equation}{0}

In this section, we recall some results established in Li (2012) on flows
of CB-processes and nonlocal branching superprocesses over the positive
half line. By a \emph{branching mechanism} $\phi$ we mean a function
$\phi$ on $[0,\infty)$ with the representation
 \beqlb\label{1.1}
\phi(z)=bz+\frac{1}{2}\sigma^2z^2 +\int_0^\infty (e^{-zu}-1+zu)m(du),
 \eeqlb
where $\sigma\geq 0$ and $b$ are constants and $(u\wedge{u}^2)m(du)$ is a
finite measure on $(0,\infty)$. Consider a family of branching mechanisms $\{\phi_q: q\in [0,1]\}$ that is
\emph{admissible} in the sense that each $\phi_q$ is given by (\ref{1.1})
with parameters $(b,m)=(b_q,m_q)$ depending on $q\in [0,1]$ and for each
$z\geq 0$ the function $q\mapsto \phi_q(z)$ is decreasing and continuously
differentiable with the derivative $\psi_\theta(z) =
-({\partial}/{\partial\theta}) \phi_{\theta}(z)$ of the form
 \beqlb\label{1.2}
\psi_\theta(z)=h_{\theta}z+\int_0^{\infty}(1-e^{-zu})n_{\theta}(du),
 \eeqlb
where $h_{\theta}\geq 0$ and $n_{\theta}(du)$ is a $\sigma$-finite kernel
from $[0,1]$ to $(0,\infty)$ satisfying
 \beqnn
\sup_{0\leq\theta\leq
1}\Big[h_{\theta}+\int_0^{\infty}un_{\theta}(du)\Big]<\infty.
 \eeqnn
Let $m(dz,d\theta)$ be the measure on $(0,\infty)\times[0,1]$ defined by
 \beqnn
m([c,d]\times [0,q])=m_q[c,d],\quad q\in [0,1], d>c>0.
 \eeqnn
Let $W(ds,du)$ be a white noise on $(0,\infty)^2$ based on the Lebesgue
measure, $\tilde{N}(ds,dz,d\theta,du)$ be a compensated Poisson random
measure on $(0,\infty)^2\times [0,1]\times (0,\infty)$ with intensity
$dsm(dz,d\theta)du$. By the results in Li (2012), the following stochastic
equation
 \beqlb\label{4.16}
Y_t(q) \ar=\ar Y_0(q)-b_q\int_0^t Y_{s-}(q)ds+
\sigma\int_0^t\int_0^{Y_{s-}(q)}W(ds,du)\cr \ar\ar
 +\int_0^t\int_0^{\infty}\int_{[0,q]}\int_0^{Y_{s-}(q)}
 z \tilde{N}(ds,dz,d\theta,du)
 \eeqlb
has a unique solution flow $\{Y_t(q):t\geq 0,q\in [0,1]\}$. For each $q\in
[0,1]$, the one-dimensional process $\{Y_t(q): t\geq 0\}$ is a CB-process
 with branching mechanism $\phi_q$. The flow is increasing in $q\in[0,1]$.
 It was verified in Li (2012) that $\{(Y_t(q))_{t\geq 0} :q\in [0,1]\}$
 can be identified as a path-valued branching process. Moreover, the flow
 induces a c\`{a}dl\`{a}g
$M[0,1]$-valued  superprocess $\{Y_t: t\ge 0\}$ which is the unique
solution of the following martingale problem: For every $G\in
C^2(\mbb{R})$ and $f\in C[0,1]$,
 \beqlb\label{4.17}
 G(\<Y_t,f\>)
 \ar=\ar
G(\<Y_0,f\>)+\, \int_0^t G^{\prime}(\<Y_s,f\>)ds
\int_{[0,1]}Y_s(dx)\int_{[0,1]}f(x\vee \theta)h_{\theta}d\theta\cr
 \ar\ar
 -\, b_0\int_0^t G^\prime(\<Y_s,f\>)\<Y_s,f\> ds
+ \frac{1}{2}\sigma^2\int_0^t G''(\<Y_s,f\>)\<Y_s,f^2\> ds \cr
 \ar\ar
+ \int_0^t ds \int_{[0,1]}Y_s(dx)\int_0^\infty \Big[G(\<Y_s,f\> + zf(x))
\cr
 \ar\ar
-\, G(\<Y_s,f\>) - zf(x)G^\prime(\<Y_s,f\>)\Big] m_0(dz) \cr
 \ar\ar
+ \int_0^t ds \int_{[0,1]}Y_s(dx)\int_{[0,1]}d\theta \int_0^\infty
\Big[G(\<Y_s,f\> + zf(x\vee \theta)) \cr
 \ar\ar
-\, G(\<Y_s,f\>)\Big] n_{\theta}(dz) +\, \mbox{local mart.}
 \eeqlb
Let $f\mapsto \Psi(\cdot,f)$ be the operator on $C^+[0,1]$ defined by
 \beqlb\label{1.4} \Psi(x,f)=\int_{[0,1]} f(x\vee \theta)h_{\theta}d \theta
+\int_{[0,1]}d\theta\int_0^{\infty}(1-e^{-zf(x\vee
\theta)})n_{\theta}(dz). \eeqlb Then the superprocess $\{Y_t: t\ge 0\}$
has local branching  mechanism $\phi_0$ and nonlocal branching  mechanism
$\Psi$. Its transition  semigroup $(Q_t)_{t\geq 0}$ is given  by
 \beqlb\label{1.5}
\int_{M[0,1]}e^{-\<\nu,f\>}Q_t(\mu,d\nu)=\exp\Big\{-\<\mu,V_tf\>\Big\},
\qquad f\in C^+[0,1], \eeqlb where $t\mapsto V_tf$ is the unique locally
bounded positive solution of \beqlb\label{1.6} V_tf(x)=f(x)-\int_0^t
[\phi_0(V_sf(x))-\Psi(x,V_sf)]ds,\qquad t\geq 0, x\in [0,1]. \eeqlb The
reader may refer Li (2012) for the derivations of the superprocess
$\{Y_t:t\geq 0\}$.


\section{Stochastic equations for discrete state branching processes}

\setcounter{equation}{0} In this section, we give a construction of the
continuous time and discrete state branching process as the solution of a
stochastic integral equation driven by Poisson random measure.
Stochastic integral equations of this type were used in Li and Ma (2008)
to construct catalytic branching processes. We here give all the details
for completeness.

Let $g=g(z)=\sum_{i=0}^\infty p_iz^i$ be a probability generating function
with $g'(1)<\infty$. Let $N(ds,dz,du)$ be a Poisson random measure on
$(0,\infty)\times \mathbb{N}\times(0,\infty)$ with intensity $\sigma
ds\pi(dz)du$, where $\sigma>0$ is a constant and
$\pi(dz):=\sum_{i=0}^\infty p_i\delta_i(dz)$. Suppose that $X_0$ is a
non-negative integer-valued random variable satisfying ${\bf
E}[X_0]<\infty$. We assume $X_0$ is independent of $N(ds,dz,du)$ and
consider the stochastic integral equation \beqlb\label{4.1}
X_t=X_0+\int_0^t\int_{\mbb{N}}\int_0^{X_{s-}}(z-1)N(ds,dz,du). \eeqlb By a
\emph{solution} of (\ref{4.1}) we mean a non-negative c\`{a}dl\`{a}g
progressive process $\{X_t:t\geq 0\}$ satisfying the equation a.s. for
each $t\geq 0$. We say \emph{pathwise uniqueness} of solution holds for
(\ref{4.1}) if any two solutions of the equation with the same initial
state are indistinguishable.

\btheorem\label{t4.1} Suppose that $\{X_t^1\}$ and $\{X_t^2\}$ are two
solutions of (\ref{4.1}) satisfying ${\bf E}[|X_0^1+X_0^2|]< \infty$. Then
we have \beqlb\label{4.2} {\bf E}[|X_t^2-X_t^1|]\leq {\bf
E}[|X_0^2-X_0^1|]\exp\{\sigma t(g'(1)+1)\}. \eeqlb Consequently, the
pathwise uniqueness of solution holds for (\ref{4.1}). \etheorem

\proof The pathwise uniqueness for (\ref{4.1}) follows from Theorem 2.1 of
Dawson and Li (2012). We present a proof of the result here for
completeness. Let $\xi_t=X_t^2-X_t^1$ for $t\geq 0$. From (\ref{4.1}) we
have \beqnn \xi_t \ar =\ar
X_0^2-X_0^1+\int_0^t\int_{\mbb{N}}\int_{X_{s-}^1}^{X_{s-}^2}
(z-1)1_{\{X_{s-}^1\leq X_{s-}^2\}}N(ds,dz,du)\cr \ar\ar
 - \int_0^t\int_{\mbb{N}}\int_{X_{s-}^2}^{X_{s-}^1}
(z-1)1_{\{X_{s-}^1> X_{s-}^2\}}N(ds,dz,du). \eeqnn Let $\tau_m=\inf\{t\geq
0: X_t^1\geq m ~or~ X_t^2\geq m\}$. Then we have \beqnn {\bf
E}[|\xi_{t\wedge\tau_m}|] \ar\leq\ar {\bf E} [|\xi_0|]+{\bf
E}\int_0^{t\wedge\tau_m}\int_{\mbb{N}}\int_{X_{s-}^1}^{X_{s-}^2}
(z+1)1_{\{X_{s-}^1\leq X_{s-}^2\}}N(ds,dz,du)\cr \ar\ar
 + {\bf E}\int_0^{t\wedge\tau_m}\int_{\mbb{N}}\int_{X_{s-}^2}^{X_{s-}^1}
(z+1)1_{\{X_{s-}^1> X_{s-}^2\}}N(ds,dz,du)\cr \ar=\ar {\bf
E}[|\xi_0|]+{\bf
E}\int_0^{t\wedge\tau_m}ds\int_{\mbb{N}}\xi_{s-}1_{\{\xi_{s-}\geq
0\}}(z+1)\sigma\pi(dz)\cr \ar\ar
 + {\bf E}\int_0^{t\wedge\tau_m}ds\int_{\mbb{N}}(-\xi_{s-})1_{\{\xi_{s-}< 0\}}(z+1)\sigma\pi(dz)\cr
 \ar\leq\ar
{\bf E}[|\xi_0|]+\int_0^{t}{\bf E}[|\xi_{s\wedge\tau_m}|]\sigma(g'(1)+1)ds.
 \eeqnn By Gronwall's inequality we get \beqnn {\bf
E}[|\xi_{t\wedge\tau_m}|]\leq {\bf E}[|\xi_0|]\exp\{\sigma t(g'(1)+1)\}.
 \eeqnn Then (\ref{4.2}) follows by Fatou's lemma.\qed

By Theorem 2.5 in Dawson and Li (2012), there is a unique strong solution
to (\ref{4.1}). Here we give a simple direct proof of the existence of the
solution. We first take an $n\in\mbb{N}_+$ and
 consider the following stochastic equation
 \beqlb\label{4.3} X_t=X_0+\int_0^t\int_{\mbb{N}}\int_0^{X_{s-}\wedge
n}(z-1)N(ds,dz,du). \eeqlb

\bproposition\label{p4.2} Let $\{X_t^n\}$ be a solution of (\ref{4.3}).
Then we have \beqlb\label{4.4} {\bf E} \Big[\sup_{0\leq s\leq
t}X_s^n\Big]\leq {\bf E}[X_0]\exp \{\sigma g'(1) t\}, \quad t\geq 0.
 \eeqlb \eproposition

\proof From (\ref{4.3}) we have \beqnn {\bf E} \Big[\sup_{0\leq s\leq
t}X_s^n\Big] \ar\leq\ar {\bf E} [X_0]+{\bf
E}\Big[\int_0^t\int_{\mbb{N}}\int_0^{X_{s-}^n\wedge n}z
N(ds,dz,du)\Big]\cr \ar=\ar
 {\bf E} [X_0]+{\bf E}\Big[\int_0^tds\int_{\mbb{N}}({X_{s-}^n\wedge n})z \sigma\pi(dz)\Big].
 \eeqnn Thus $t\mapsto {\bf E}[\sup_{0\leq s\leq t}X_s^n]$ is a locally
bounded function. Moreover, \beqnn {\bf E} \Big[\sup_{0\leq s\leq
t}X_s^n\Big] \ar\leq\ar
 {\bf E} [X_0]+\int_0^t ds\int_{\mbb{N}}{\bf E}\Big[\sup_{0\leq r\leq s} X_r^n\Big]z \sigma\pi(dz)\cr
 \ar=\ar
 {\bf E} [X_0]+\sigma g'(1)\int_0^t  {\bf E}\Big[\sup_{0\leq r\leq s} X_r^n\Big] ds.
 \eeqnn By Gronwall's lemma we get the result.\qed

By a modification of the proof of Theorem \ref{t4.1} we get the following Proposition.
\bproposition\label{t4.3} Suppose that $\{X_t^{n,1}\}$ and $\{X_t^{n,2}\}$
are two solutions of (\ref{4.3}). Then we have \beqlb\label{4.5} {\bf
E}[|X_t^{n,2}-X_t^{n,1}|]\leq {\bf E}[|X_0^{n,2}-X_0^{n,1}|]\exp\{\sigma
t(g'(1)+1)\}. \eeqlb Consequently, the pathwise uniqueness of solution
holds for (\ref{4.3}). \eproposition

\bproposition\label{p4.4} For each $n\geq 1$, there is a solution
$\{X_{t}^n:t\geq 0\}$ of (\ref{4.3}). \eproposition

\proof Let $\{S_k:k=1,2,\cdots\}$ be the set of jump times of the Poisson
process \beqnn t\mapsto \int_0^t\int_{\mbb{N}}\int_0^n N(ds,dz,du).
 \eeqnn We have clearly $S_k\to\infty$ as $k\to\infty$. For $0\leq t<S_1$,
set $X_t^n=X_0$. Suppose that $X_t^n$ has been defined for $0\leq t<S_k$
and let \beqnn
X_t^n=X_{S_{k}-}^n+\int_{\{S_k\}}\int_{\mbb{N}}\int_0^{X_{S_{k}-}^n\wedge
n} (z-1)N(ds,dz,du), \quad S_k\leq t<S_{k+1}. \eeqnn From the construction
of $X_{S_k}^n$ we see $X_{S_k}^n-X_{S_{k-1}}^n\geq -1$.  And since
$X_{S_{k-1}}^n=0$ implies $X_{S_k}^n=0$, $X_{S_k}^n\in \mbb{N}$. By
induction that defines a non-negative process $\{X_{t}^n:t\geq 0\}$ which
is clearly a solution to (\ref{4.3}).\qed

\bproposition\label{p4.5} Let $\{X_{t}^n: t\geq 0\}$ be the solution of
(\ref{4.3}) with $n=1,2,\cdots$. Then the sequence $\{X_{t}^n:t\geq 0\}$
is tight in $D([0,\infty),\mbb{N})$. \eproposition

\proof By Proposition \ref{p4.2}, it is easy to see that \beqnn t\mapsto
C_t:=\sup_{n\geq 1} {\bf E}\Big[\sup_{0\leq s\leq t}X_s^n\Big] \eeqnn is
locally bounded. Then for every fixed $t\geq 0$, the sequence of random
variables $X_{t}^n$ is tight. Moreover, in view of (\ref{4.3}), if
$\{\tau_n\}$ is a sequence of stopping times bounded above by $T\geq 0$, we
have \beqnn {\bf E}[|X_{t+\tau_n}^{n}-X_{\tau_n}^{n}|] \ar=\ar {\bf
E}\Big[\int_{\tau_n}^{t+\tau_n} \int_{\mbb{N}}\int_0^{{X_{s-}^n}\wedge n}
(z-1) N(ds,dz,du)\Big]\cr \ar\leq\ar {\bf E}\Big[\int_{0}^{t}
ds\int_{\mbb{N}}({X_{s+\tau_n}^n\wedge n})(z+1) \sigma \pi(dz)\Big]\cr
\ar\leq\ar \sigma(g'(1)+1)\int_0^t {\bf E}[X_{s+\tau_n}^n]ds\cr \ar\leq\ar
{\bf E}[X_0]\exp\{\sigma g'(1)(t+T)\}\sigma(g'(1)+1)t, \eeqnn where the
last inequality follows by Proposition \ref{p4.2}. Consequently, as $t\to
0$, \beqnn \sup_{n\geq 1} {\bf E}[|X_{t+\tau_n}^{n}-X_{\tau_n}^{n}|]\to 0.
\eeqnn Then $\{X_{t}^n:t\geq 0\}$ is tight in $D([0,\infty),\mbb{N})$ by
the criterion of Aldous (1978); see also Ethier and Kurtz (1986,
pp.137-138). \qed \btheorem\label{t4.6} There is a solution $\{X_{t}:t\geq
0\}$ of (\ref{4.1}). \etheorem

\proof For each $n\geq 1$, let $\{X_{t}^n:t\geq 0\}$ be the solution of
(\ref{4.3}). Define $\tau_n=\inf\{t\geq 0:X_t^n\geq n\}$. From Proposition
\ref{p4.2} it follows that \beqnn {\bf E} [X_{t\wedge\tau_n}^n]\leq {\bf
E} \Big[\sup_{0\leq s\leq t}X_s^n\Big] \leq {\bf E}[X_0]\exp \{\sigma
g'(1) t\}, \quad t\geq 0. \eeqnn Then we have \beqnn {\bf E}
[X_{t\wedge\tau_n}^n1_{\{\tau_n\leq t\}}]\leq {\bf E}[X_0]\exp \{\sigma
g'(1) t\}. \eeqnn By the right continuity of $\{X_t^n\}$ we have
$X_{\tau_n}^n\geq n$, so \beqnn n{\bf P}[\{\tau_n\leq t\}]\leq {\bf
E}[X_0]\exp \{\sigma g'(1) t\},  \quad t\geq 0. \eeqnn That implies
$\tau_n\to\infty$ almost surely as $n\to\infty$. On the other hand,
$\{X_{t}^n\}$ satisfies the equation (\ref{4.1}) for $0\leq t<\tau_n$. By
the pathwise  uniqueness of the solution of  (\ref{4.1}) we get, for any
$i,j\in\mbb{N}$,
 \beqnn
 X_t^i=X_t^j, \quad t<\tau_i\wedge\tau_j.
 \eeqnn
Let $\{X_{t}\}$ be the process such that $X_t=X_t^n$ for all $0\leq
t<\tau_n$ and $n\geq 1$. It is easily seen that $\{X_{t}\}$ is a solution
of (\ref{4.1}).\qed

 Theorems \ref{t4.1} and \ref{t4.6} imlpy that (\ref{4.1}) has unique \emph{strong solution} and the
 solution $\{X_{t}:t\geq 0\}$ is a strong Markov process; see, e.g., Ikeda and Watanabe (1989, pp.163-166 and p.215).
Let $B(\mbb{N})$ denote the set of bounded measurable functions on
$\mbb{N}$.
 By It\^{o}'s formula it is easy to see that  $\{X_{t}:t\geq 0\}$ has generator $A$ defined by
 \beqnn
 Af(x)=\sigma x\sum_{i=0}^{\infty}[f(x+i-1)-f(x)]p_i,\qquad x\in\mbb{N}, ~f\in B(\mbb{N}).
 \eeqnn
 Then $\{X_{t}:t\geq 0\}$ is a continuous time Galton-Watson branching process.

 Let $N^{(1)}(ds,dz,du)$ and $N^{(2)}(ds,dz,du)$ be two mutually independent Poisson random measures on
 $(0,\infty)\times \mathbb{N}\times(0,\infty)$ with the same intensity $\sigma ds\pi(dz)du$.
 Consider the following two stochastic equations
 \beqnn
X_t^{(1)}=X_0^{(1)}+\int_0^t\int_{\mbb{N}}\int_0^{{X_{s-}^{(1)}}}(z-1)N^{(1)}(ds,dz,du)
\eeqnn
and
\beqnn
X_t^{(2)}=X_0^{(2)}+\int_0^t\int_{\mbb{N}}\int_0^{{X_{s-}^{(2)}}}(z-1)N^{(2)}(ds,dz,du).
 \eeqnn
Clearly,  $X_t^{(1)}$ and  $X_t^{(2)}$ are mutually independent. Set
$X_t=X_t^{(1)}+X_t^{(2)}$. Since the random measure
 \beqnn
 N'(ds,dz):=
\int_{\{0<u\leq {X_{s-}^{(1)}}\}} N^{(1)}(ds,dz,du)+\int_{\{0<u\leq
{X_{s-}^{(2)}}\}} N^{(2)}(ds,dz,du)
 \eeqnn
has predictable compensator $\sigma X_{s-} ds\pi(dz)$,  by representation
theorems for semimartingales, on an extension of the original probability
space, there is a Poisson random measure on
 $(0,\infty)\times \mathbb{N}\times(0,\infty)$ with  intensity $\sigma ds\pi(dz)du$ such that
 \beqnn X_t=X_0+\int_0^t\int_{\mbb{N}}\int_0^{X_{s-}}(z-1)N(ds,dz,du);
 \eeqnn see, e.g., Ikeda and Watanabe (1989, p.93). Then the solution of
(\ref{4.1}) is a branching process (continuous time and discrete  state).
This gives another derivation of the branching property of $\{X_{t}:t\geq
0\}$.


\section{The flow of discrete state branching processes}

\setcounter{equation}{0} In this section, we give a formulation of the discrete state
  branching flow as the solution flow of a set of stochastic integral equations.
Let $\{g_{\theta}:\theta\geq 0\}$ be a family of  probability generating
functions, that is, for each $\theta\geq 0$, \beqnn
g_{\theta}(z)=\sum_{i=0}^\infty p_i(\theta)z^i,\quad |z|\leq 1, \eeqnn
where $p_i(\theta)\geq 0$ and $\sum_{k=0}^{\infty}p_i(\theta)=1$.
Moreover, we assume $\theta\mapsto g_\theta'(1)$ is continuous and
$p_i({\theta_2})\geq p_i({\theta_1})$
 holds for any $\theta_2\geq \theta_1\geq 0$ and $i\in \mbb{N}_+$.
Define a family of probability measures $\{\pi_{\theta}:\theta\geq 0\}$ on
$\mbb{N}$ by \beqnn \pi_{\theta}(dz)=\sum_{i=0}^\infty
p_i(\theta)\delta_i(dz) \eeqnn Then we have
$\pi_{\theta_2}|_{\mbb{N_+}}\geq \pi_{\theta_1}|_{\mbb{N_+}}$ for any
$\theta_2\geq \theta_1\geq 0$. Let $\bar{\pi}(dz,d\theta)$ be the measure
on $\mbb{N}_+\times [0,\infty)$ defined by \beqnn \bar{\pi}(A\times
[0,\theta])=\pi_{\theta}(A),\quad A\subset \mbb{N}_+, \theta\geq 0. \eeqnn
Notice that the positive function $\theta\mapsto
b(\theta):=\pi_{\theta}(\{0\})$ is decreasing.

Let $q\mapsto X_0(q)$ be a deterministic positive right continuous
increasing function on $[0,\infty)$ and take values in $\mbb{N}$. Let
$N(ds,dz,d\theta,du)$ be a Poisson random measure on
$(0,\infty)\times{\mbb{N}}\times[0,\infty)\times(0,\infty)$ with intensity
$\sigma ds\bar{\pi}(dz,d\theta)du$ and $N_0(ds,d\theta,du)$ a Poisson
random measure on $(0,\infty)^3$ with intensity $\sigma ds\d\theta du$.
Suppose that $N(ds,dz,d\theta,du)$ and $N_0(ds,d\theta,du)$ are
independent of each other. Consider stochastic integral equation
 \beqlb\label{4.6} X_t(q) \ar=\ar
X_0(q)+\int_0^t\int_{\mbb{N}_+}\int_{[0,q]}\int_0^{X_{s-}(q)}(z-1)N(ds,dz,d\theta,du)\cr
\ar\ar -\int_0^t\int_{0}^{b(q)}\int_0^{X_{s-}(q)}N_0(ds,d\theta,du).
 \eeqlb Note that for each $q\geq 0$, \beqnn \int_{\{0<\theta\leq b(q)\}}
N_0(ds,d\theta,du) \eeqnn is a Poisson random measure with intensity
$\sigma b(q)dsdu=\sigma \bar{\pi}_0(\mbb{N}\times[0,q])dsdu$, where
$\bar{\pi}_0(dz,d\theta)$ is a measure on $\mbb{N}\times[0,\infty)$
defined by \beqnn
\bar{\pi}_0(A\times[0,q])=\pi_{q}(\{0\})\delta_0(A),\quad A\subset\mbb{N},
~\theta\geq 0. \eeqnn By representation theorems for semi-martingales,
there is a Poisson random measure $N_1(ds,dz,d\theta,du)$ on
$(0,\infty)\times{\mbb{N}}\times[0,\infty)\times(0,\infty)$ with intensity
$\sigma ds\bar{\pi}_0(dz,d\theta)du$ such that for every
$E\in\mathscr{B}(0,\infty)$, \beqnn \int_0^t\int_0^{b(q)}\int_E
N_0(ds,d\theta,du)= \int_0^t\int_{\mbb{N}}\int_{[0,q]}\int_E
N_1(ds,dz,d\theta,du); \eeqnn
 see, e.g., Ikeda and Watanabe (1989, p.93). Define $N_2(ds,dz,du)$ by
 \beqnn
 N_2(ds,dz,du)=\int_{\{0\leq \theta\leq q\}}N(ds,dz,d\theta,du)+\int_{\{0\leq \theta\leq q\}}N_1(ds,dz,d\theta,du).
 \eeqnn
  Then $N_2$ is a Poisson random measure on $(0,\infty)\times{\mbb{N}}\times(0,\infty)$
  with intensity $\sigma ds\pi_q(dz)du$ and the equation (\ref{4.6}) can be rewrited as
 \beqnn
  X_t(q)=X_0(q)+\int_0^t\int_{\mbb{N}}\int_0^{X_{s-}(q)}(z-1) N_2(ds,d\theta,du).
 \eeqnn By Theorem \ref{t4.6} we see that for each $q\geq 0$, the equation
(\ref{4.6}) has a unique strong solution $\{X_t(q):t\geq 0\}$.

 \btheorem\label{t4.7}
Suppose that $q\geq p\geq 0$. Let $\{X_t(q)\}$ be the solution of
(\ref{4.6}) and $\{X_t(p)\}$ be the solution of the equation with $q$
replaced by $p$. Then we have ${\bf P}\{X_t(q)\geq X_t(p)~ for ~all ~t\geq
0\}=1$. \etheorem

 \proof
 Let $\zeta_t=X_t(p)-X_t(q)$ for $t\geq 0$. From (\ref{4.6}) we have
 \beqlb\label{4.7}
 \zeta_t
 \ar=\ar
 \zeta_0+\int_0^t\int_{\mbb{N}_+}\int_{[0,p]}\int_{X_{s-}(q)}^{X_{s-}(p)}(z-1)N(ds,dz,d\theta,du)\cr
\ar\ar
-\int_0^t\int_{\mbb{N}_+}\int_{(p,q]}\int_0^{X_{s-}(q)}(z-1)N(ds,dz,d\theta,du)
-\int_0^t\int_{0}^{b(q)}\int_{X_{s-}(q)}^{X_{s-}(p)}N_0(ds,d\theta,du)\cr
\ar\ar -\int_0^t\int_{b(q)}^{b(p)}\int_0^{X_{s-}(p)}N_0(ds,d\theta,du).
 \eeqlb
 Let $\tau_m=\inf\{t\geq 0: X_t(q)\geq m ~or~ X_t(p)\geq m\}$.
It is easy to construct a sequence of functions $\{f_n\}$ on $\mbb{R}$
such that $0\leq f_n'(z)\leq 1$ for $z\geq 0$ and $f_n(z)=f_n'(z)=0$ for
$z\leq 0$. Moreover, $f_n(z)\to z^+:=0\vee z$ increasingly as
$n\to\infty$. By (\ref{4.7}) and It\^{o}'s formula, \beqnn
f_n(\zeta_{t\wedge \tau_m}) \ar=\ar \int_0^{t\wedge
\tau_m}\int_{\mbb{N}_+}\int_{[0,p]}\int_{X_{s-}(q)}^{X_{s-}(p)}
[f_n(\zeta_{s-}+z-1)-f_n(\zeta_{s-})]1_{\{\zeta_{s-}>
0\}}N(ds,dz,d\theta,du)\cr \ar\ar +\int_0^{t\wedge
\tau_m}\int_{\mbb{N}_+}\int_{(p,q]}\int_0^{X_{s-}(q)}
[f_n(\zeta_{s-}-z+1)-f_n(\zeta_{s-})]N(ds,dz,d\theta,du)\cr \ar\ar
+\int_0^{t\wedge \tau_m}\int_0^{b(q)}\int_{X_{s-}(p)}^{X_{s-}(q)}
[f_n(\zeta_{s-}-1)-f_n(\zeta_{s-})]1_{\{\zeta_{s-}>
0\}}N_0(ds,d\theta,du)\cr \ar\ar +\int_0^{t\wedge
\tau_m}\int_{b(q)}^{b(p)}\int_0^{X_{s-}(p)}
[f_n(\zeta_{s-}-1)-f_n(\zeta_{s-})]N_0(ds,d\theta,du)\cr \ar\leq\ar \sigma
\int_0^{t\wedge \tau_m}\zeta_{s-}1_{\{\zeta_{s-}>
0\}}ds\int_{\mbb{N}_+}(z-1) \pi_p(dz) +\mbox{martingale}. \eeqnn
 Taking the expectation in both sides and letting $n\to\infty$ gives
 \beqnn
 {\bf E}[\zeta_{t\wedge\tau_m}^+]\leq \sigma (g_p'(1)-1+b(p))\int_0^t {\bf E}[\zeta_{s\wedge\tau_m}^+]ds.
 \eeqnn
 Then ${\bf E}[\zeta_{t\wedge\tau_m}^+]=0$ for all $t\geq 0$. Since $\tau_m\to\infty$ as $m\to\infty$,
 that proves the desired comparison result.\qed
\bproposition\label{p4.8}
There is a locally bounded  positive function $(t,u)\mapsto C(t,u)$ on
$[0,\infty)^2$ so that, for any $t\geq 0$ and $p\leq q\leq u<\infty$,
 \beqlb\label{4.8} {\bf E} \Big\{\sup_{0\leq s\leq
t}[X_s(q)-X_s(p)]\Big\}\leq
C(t,u)\Big\{X_0(q)-X_0(p)+g_q'(1)-g_p'(1)\Big\}. \eeqlb \eproposition

 \proof
 Let $\xi_t=X_t(q)-X_t(p)$. From (\ref{4.6}) we get
 \beqnn
 \sup_{0\leq s\leq t}\xi_s
 \ar\leq\ar
 \xi_0+\int_0^t\int_{\mbb{N}_+}\int_{[0,q]}\int_{X_{s-}(p)}^{X_{s-}(q)}(z-1)N(ds,dz,d\theta,du)\cr
\ar\ar
+\int_0^t\int_{\mbb{N}_+}\int_{(p,q]}\int_0^{X_{s-}(p)}(z-1)N(ds,dz,d\theta,du)\cr
\ar\ar
+\int_0^t\int_{b(q)}^{b(p)}\int_0^{X_{s-}(p)}N_0(ds,d\theta,du).
 \eeqnn
 Then
 \beqnn
 {\bf E}\Big[\sup_{0\leq s\leq t}\xi_s\Big]
 \ar\leq\ar
\xi_0+\sigma [g_p'(1)-1+b(q)]\int_0^t {\bf E}[\xi_s]ds\cr \ar\ar
+\sigma[g_q'(1)-g_p'(1)]\int_0^t{\bf E}[X_s(p)]ds.
 \eeqnn
 Since $t\mapsto {\bf E}[X_t(p)]$ is locally bounded, by Gronwall's inequality we get the desired estimate.\qed

 From the discussion above, given a constant $\sigma>0$ and
a family of  probability generating functions $\{g_{\theta}:\theta\geq
0\}$, we obtain a  continuous time and discrete state branching process
flow
 $\{X_t(q):t\geq 0,q\geq 0\}$ as the solution of equation (\ref{4.6}).
For any $t\geq 0$ define the random function
 $\tilde{X}_t\in F[0,1]$ by
 $\tilde{X}_t(1)={X}_t(1)$ and
 \beqlb\label{4.19}
 \tilde{X}_t(q)=\inf\{{X}_t(u):\mbox{rational} ~u\in(q,1]\},\quad 0\leq q<1.
 \eeqlb
 By Proposition \ref{p4.8}, for each $q\in[0,1]$ we have
 \beqnn
 {\bf P}\{\tilde{X}_t(q)=X_t(q) ~\mbox{for all}~ t\geq 0\}=1.
 \eeqnn
 Then $\{\tilde{X}_t(q):t\geq 0\}$ is also c\`{a}dl\`{a}g and solves (\ref{4.6}) for every $q\in[0,1]$.


\section{Scaling limits of the discrete branching flows}

\setcounter{equation}{0}

In this section, we prove some limit theorems for the discrete state
branching flows, which will lead to the continuous state branching flows
of Li (2012). We shall present the limit theorems in the settings of
measure-valued processes and path-valued processes.

Suppose that for each $k\geq 1$, there is a positive constant $\sigma_k$
and a family of generating functions $\{g_{\theta}^{(k)}: \theta\geq 0\}$
satisfying the assumptions specified at the beginning of the last section.
Then we can define $\pi_{\theta}^{(k)}(dz)$ and $\bar{\pi}^{(k)}
(dz,d\theta)$ in the same way as there. Moreover, assume $\theta\mapsto
b_k(\theta) := g_{\theta}^{(k)}(0)$ is differentiable. Let
$\{X_t^{(k)}(q):t\geq 0\}$ be the corresponding solution of (\ref{4.6})
and $\{\tilde{X}_t^{(k)}(q):t\geq 0, q\in[0,k]\}$ be defined in the same
way as in (\ref{4.19}). Define
\beqlb\label{4.9a}
Y_t^{(k)}(q)=\frac{1}{k}\tilde{X}_{t}^{(k)}(kq),\qquad q\in [0,1].
\eeqlb
From (\ref{4.6}) we have
 \beqlb\label{4.9}
 Y_t^{(k)}(q)
  \ar=\ar
Y_0^{(k)}(q)+\frac{1}{k}\int_0^t\int_{\mbb{N}_+} \int_{[0,kq]}
\int_0^{kY_{s-}^{(k)}(q)} (z-1) N(ds,dz,d\theta,du)\cr
 \ar\ar
-\frac{1}{k}\int_0^t\int_{0}^{b_k(kq)}\int_0^{kY_{s-}^{(k)}(q)}N_0(ds,d\theta,du).
 \eeqlb One can use a standard stopping time argument to show that for any
$q\in [0,1]$,
 the function $t\mapsto \mbf{E}[Y_t^{(k)}(q)]$ is locally bounded.
Then by an  argument similar to the proof of Proposition \ref{p4.2} we
have
 \bproposition\label{p4.9}
 For any $t\geq 0$ and $q\in [0,1]$, we have
 \beqlb\label{4.10} {\bf E} \Big[\sup_{0\leq s\leq t}Y_s^{(k)}(q)\Big] \leq
Y_0^{(k)}(q)\exp\Big\{t\sigma_k
\Big((g_{kq}^{(k)})'(1)-1+b_k(kq)\Big)\Big\}. \eeqlb \eproposition

The random function $Y_t^{(k)}\in F[0,1]$ induces a random measure
$Y_t^{(k)}\in M[0,1]$ so that
 $Y_t^{(k)}([0,q])=Y_t^{(k)}(q)$ for $q\in [0,1]$. We are interested in the asymptotic behavior of
 $\{Y_t^{(k)}:t\geq 0\}$ as $k\to\infty$. For any $f\in C^1[0,1]$
 one can use Fubini's theorem to see
 \beqlb\label{4.11}
 \<Y_t^{(k)},f\>=f(1)Y_t^{(k)}(1)-\int_0^1 f'(q)Y_t^{(k)}(q)dq.
 \eeqlb
Fix an integer $n\geq 1$ and let $q_i=i/2^n$ for $i=0,1,\cdots,2^n$. By
(\ref{4.9}) we have \beqlb\label{4.12}
\sum_{i=1}^{2^n}f'(q_i)Y_t^{(k)}(q_i) \ar=\ar
\sum_{i=1}^{2^n}f'(q_i)Y_0^{(k)}(q_i)\cr \ar\ar
+\frac{1}{k}\sum_{i=1}^{2^n}f(q_i)\int_0^t\int_{\mbb{N}_+}\int_{[0,kq_i]}\int_0^{kY_{s-}^{(k)}(q_i)}(z-1)N(ds,dz,d\theta,du)\cr
\ar\ar
-\frac{1}{k}\sum_{i=1}^{2^n}f'(q_i)\int_0^t\int_{0}^{b_k(kq_i)}\int_0^{kY_{s-}^{(k)}(q_i)}N_0(ds,d\theta,du)\cr
\ar=\ar \sum_{i=1}^{2^n}f'(q_i)Y_0^{(k)}(q_i)\cr \ar\ar
+\frac{1}{k}\int_0^t\int_{\mbb{N}_+}\int_{[0,k]}\int_0^{kY_{s-}^{(k)}(1)}F_n^{(k)}(s,\theta,u)(z-1)N(ds,dz,d\theta,du)\cr
\ar\ar
-\frac{1}{k}\int_0^t\int_{0}^{b_k(0)}\int_0^{kY_{s-}^{(k)}(1)}\tilde{F}_n^{(k)}(s,\theta,u)N_0(ds,d\theta,du),
 \eeqlb
 where
 \beqnn
 F_n^{(k)}(s,\theta,u)=\sum_{i=1}^{2^n}f'(q_i)1_{\{\theta\leq kq_i\}}1_{\{u\leq kY_{s-}^{(k)}(q_i)\}}
 \eeqnn
 and
 \beqnn
\tilde{F}_n^{(k)}(s,\theta,u)=\sum_{i=1}^{2^n}f'(q_i)1_{\{\theta\leq
b_k(kq_i)\}}1_{\{u\leq kY_{s-}^{(k)}(q_i)\}}.
 \eeqnn
By the right continuity of $q\mapsto Y_t^{(k)}(q)$ it is easy to see that,
as $n\to\infty$, \beqnn 2^{-n}F_n^{(k)}(s,\theta,u)\to
F^{(k)}(s,\theta,u) :=\int_0^1 f'(q)1_{\{\theta\leq kq\}}1_{\{u\leq
kY_{s-}^{(k)}(q)\}}dq \eeqnn and \beqnn
2^{-n}\tilde{F}_n^{(k)}(s,\theta,u)\to \tilde{F}^{(k)}(s,\theta,u)
:=\int_0^1 f'(q)1_{\{\theta\leq b_k(kq)\}}1_{\{u\leq
kY_{s-}^{(k)}(q)\}}dq. \eeqnn Then by (\ref{4.12}) we have, almost surely,
 \beqlb\label{4.13} \int_0^1 f'(q)Y_t^{(k)}(q)dq \ar=\ar \int_0^1
f'(q)Y_0^{(k)}(q)dq\cr \ar\ar
+\frac{1}{k}\int_0^t\int_{\mbb{N}_+}\int_{[0,k]}\int_0^{kY_{s-}^{(k)}(1)}
F^{(k)}(s,\theta,u)(z-1)N(ds,dz,d\theta,du)\cr \ar\ar
-\frac{1}{k}\int_0^t\int_{0}^{b_k(0)}\int_0^{kY_{s-}^{(k)}(1)}\tilde{F}^{(k)}(s,\theta,u)N_0(ds,d\theta,du).
 \eeqlb From (\ref{4.9}), (\ref{4.11}) and (\ref{4.13}) it follows that,
almost surely, \beqlb\label{4.14} \<Y_t^{(k)},f\> \ar=\ar
\<Y_0^{(k)},f\>\cr \ar\ar +
\frac{1}{k}\int_0^t\int_{\mbb{N}_+}\int_{[0,k]}\int_0^{kY_{s-}^{(k)}(1)}
[f(1)-F^{(k)}(s,\theta,u)](z-1)N(ds,dz,d\theta,du)\cr \ar\ar
-\frac{1}{k}\int_0^t\int_{0}^{b_k(k)}\int_0^{kY_{s-}^{(k)}(1)}[f(1)-\tilde{F}^{(k)}(s,\theta,u)]N_0(ds,d\theta,du)\cr
\ar\ar
+\frac{1}{k}\int_0^t\int_{{b_k(k)}}^{b_k(0)}\int_0^{kY_{s-}^{(k)}(1)}\tilde{F}^{(k)}(s,\theta,u)N_0(ds,d\theta,du).
 \eeqlb \bproposition\label{p4.10} Suppose that $Y_0^{(k)}(1)$ converges to
some $Y_0(1)$ as $k\to\infty$ and \beqnn \sup_{k\geq 1}\sigma_k
\Big[(g_{k}^{(k)})'(1)-1+b_k(0)\Big]< \infty. \eeqnn Then
$\{Y_t^{(k)}:t\geq 0\}$, $k=1,2,\cdots$ is a tight sequence in
$D([0,\infty),M[0,1])$. \eproposition

\proof For any $t\geq 0$ and $f\in C[0,1]$, by Proposition {\ref{p4.9}} it
is easy to see that \beqnn t\mapsto C_t:=\sup_{k\geq 1}{\bf E}
\Big[\sup_{0\leq s\leq t}\<Y_s^{(k)},f\>\Big] \eeqnn is locally bounded.
Then for every fixed $t\geq 0$, the sequence $\<Y_t^{(k)},f\>$ is tight.
Let $\tau_k$ be a bounded stopping time for $\{Y_t^{(k)}: t\ge 0\}$ and
assume the sequence $\{\tau_k: k=1,2,\cdots\}$ is bounded above by $T\geq
0$. Let $f\in C^1[0,1]$. By (\ref{4.14})  we see
 \beqlb\label{9.9}
\ar\ar\mbf{E}\Big[\Big|\<Y_{\tau_k+t}^{(k)},f\>
-\<Y_{\tau_k}^{(k)},f\>\Big|\Big] \cr \ar\ar\quad \leq\frac{\sigma_k }{k}
{\bf E}
\bigg[\int_0^tds\int_{\mbb{N}_+}\int_{[0,k]}\int_0^{kY_{s+\tau_k}^{(k)}(1)}
(z-1)|f(1)-F^{(k)}(s+\tau_k,\theta,u)|\bar{\pi}^{(k)}(dz,d\theta)du\bigg]\cr
\ar\ar\qquad
 +\frac{\sigma_k }{k}{\bf E}\bigg[\int_0^t ds\int_{0}^{b_k(k)}d\theta\int_0^{kY_{s+\tau_k}^{(k)}(1)}
 |f(1)-\tilde{F}^{(k)}(s+\tau_k,\theta,u)|  du\bigg]\cr
 \ar\ar\qquad
 +\frac{\sigma_k }{k}{\bf E}\bigg[\int_0^t ds\int_{{b_k(k)}}^{b_k(0)}d\theta\int_0^{kY_{s+\tau_k}^{(k)}(1)}
 |\tilde{F}^{(k)}(s+\tau_k,\theta,u)| du\bigg].
 \eeqlb
 For $s,\theta,u>0$ let $Y_{s,k}^{-1}(u)=\inf\{q\geq 0:Y_s^{(k)}(q)>u\}$
 and $b_k^{-1}(u)=\inf\{q\geq 0:b_k(q)>u\}$. It is easy to see that
 $\{q\geq 0:u\leq kY_s^{(k)}(q)\}=[Y_{s,k}^{-1}(u/k),\infty)$ and
 $\{q\geq 0:\theta\leq b_k(kq)\}=[0,b_k^{-1}(\theta)/k]$ except for at most countably many
 $u>0$ and $\theta>0$, respectively. Then in the above we can replace
 $f(1)-F^{(k)}(s,\theta,u)$ by
 \beqnn
 f(1)-\int_{\theta/k}^1 f'(q)1_{\{Y_{s,k}^{-1}(u/k)\leq q\}}dq
 =f\Big(Y_{s,k}^{-1}(\frac{u}{k})\vee \frac{\theta}{k}\Big)
 \eeqnn
 and  $\tilde{F}^{(k)}(s,\theta,u)$ can be replaced by
 \beqnn
 \ar\ar\int_0^1f'(q)1_{\{q\leq b_k^{-1}(\theta)/k\}}1_{\{Y_{s,k}^{-1}(u/k)\leq q\}}dq\cr
\ar\ar\qquad
 =\Big[f(1\wedge( b_k^{-1}(\theta)/k))-f(Y_{s,k}^{-1}(u/k))\Big]
 1_{\{Y_{s,k}^{-1}(u/k)\leq  b_k^{-1}(\theta)/k\}}.
 \eeqnn
Then from (\ref{9.9}) we have
 \beqlb\label{4.15}
\ar\ar\mbf{E}\Big[\Big|\<Y_{\tau_k+t}^{(k)},f\>
-\<Y_{\tau_k}^{(k)},f\>\Big|\Big]\cr \ar\ar\quad \leq\sigma_k {\bf E}
\bigg[\int_0^t ds\int_0^1
Y_{s+\tau_k}^{(k)}(dx)\int_{\mbb{N}_+}\int_{[0,k]} (z-1)|f(x\vee
\theta)|\bar{\pi}^{(k)}(dz,d\theta)\bigg]\cr \ar\ar\qquad
 +\sigma_k{\bf E}\bigg[\int_0^t ds\int_{0}^{b_k(k)}d\theta\int_0^1
 |f(x)|  Y_{s+\tau_k}^{(k)}(dx)\bigg]\cr
 \ar\ar\qquad
 +\sigma_k{\bf E}\bigg[\int_0^t ds\int_{{b_k(k)}}^{b_k(0)}d\theta\int_0^1
 |f(b_k^{-1}(\theta)/k)-f(x)| Y_{s+\tau_k}^{(k)}(dx)\bigg]\cr
 \ar\ar\quad
\leq \|f\|\sigma_k \int_0^t {\bf E} \Big[Y_{s+\tau_k}^{(k)}(1)\Big] ds
\int_{\mbb{N}_+}(z-1)\pi_{k}^{(k)}(dz)\cr \ar\ar\qquad
 + \|f\| \sigma_k  {b_k(k)}{\bf E}\int_0^t \Big[Y_{s+\tau_k}^{(k)}(1)\Big]ds\cr
 \ar\ar\qquad
 +2\|f\|\sigma_k[b_k(0)-b_k(k)]\int_0^t {\bf E}\Big[Y_{s+\tau_k}^{(k)}(1)\Big] ds\cr
 \ar\ar\quad
\leq \|f\| \sigma_k\Big((g_{k}^{(k)})'(1)-1+2b_k(0)\Big) \int_0^t  {\bf
E} \Big[Y_{s+\tau_k}^{(k)}(1)\Big] ds
 \cr
 \ar\ar\quad
\leq 2\|f\| Y_0^{(k)}(1)t\sigma_kA_k\exp\Big\{\sigma_k A_k(t+T)\Big\},
 \eeqlb where $A_k=(g_{k}^{(k)})'(1)-1+b_k(0)$ and the last inequality
follows by Proposition {\ref{p4.9}}. For $f\in C[0,1]$ the above
inequality follows by an approximation argument. Then we have \beqnn
\lim_{t\to 0}\sup_{k\geq 1}\mbf{E}\Big[\Big|\<Y_{\tau_k+t}^{(k)},f\>
 -\<Y_{\tau_k}^{(k)},f\>\Big|\Big]=0.
 \eeqnn By a criterion of Aldous (1978), the sequence $\{\<Y_t^{(k)},f\>:
t\geq 0\}$ is tight in $D([0,\infty),\mbb{R})$; see also Ethier and Kurtz
(1986, pp.137-138). Then the tightness criterion of Roelly (1986) implies
$\{Y_t^{(k)}: t\ge 0\}$ is tight in $D([0,\infty), M[0,1])$. \qed

For any $z\geq 0$ define
 \beqlb\label{4.4c}
\phi_{\theta}^{(k)}(z) = k\sigma_k\Big[g_{k\theta}^{(k)}(e^{-z/k})-e^{-z/k}\Big].
 \eeqlb
Let us consider the following condition:

\noindent\emph{{\bf Condition (5.A)} ~For each $l\geq 0$ the sequence $\{\phi_{\theta}^{(k)}(z)\}$ is
Lipschitz with respect to $z$ uniformly on $[0,l]^2$ and
there is an admissible family of branching mechanisms $\{\phi_{\theta}^{(k)}(z): \theta\geq 0\}$
with $(\partial/\partial\theta) \phi_{\theta}(z)=- \psi_{\theta}(z)$
such that $\phi_{\theta}^{(k)}(z)\to\phi_{\theta}(z)$  uniformly on $[0,l]^2$ as $k\to \infty$.}

Let $\{Y_t: t\ge 0\}$ be the c\`{a}dl\`{a}g superprocess with transition semigroup
defined by (\ref{1.5}) and (\ref{1.6}).

\btheorem\label{t4.11}
Suppose that Condition (5.A) holds and $\sup_{k\geq 1}\sigma_k b_k(0)<\infty$.
If $Y_0^{(k)}$ converges weakly to $Y_0\in M[0,1]$, then
$\{Y_t^{(k)}: t\ge 0\}$ converges to the superprocess
$\{Y_t: t\ge 0\}$ in distribution on $D([0,\infty), M[0,1])$.
\etheorem

 \proof
 Under the assumption, we  have
 \beqnn\sup_{k\geq 1}\sigma_k
\Big[(g_{k}^{(k)})'(1)-1+b_k(0)\Big]< \infty. \eeqnn
 By Proposition \ref{p4.10} and Skorokhod's representation theorem,  to simplify the notation we pass
to a subsequence and simply assume  $\{Y_t^{(k)}: t\ge 0\}$
converges a.s.\ to a process $\{Z_t: t\ge 0\}$ in the topology of $D([0,\infty), M[0,1])$.
 Since the solution of the martingale problem
(\ref{4.17}) is unique, it suffices to prove the weak limit point
$\{Z_t: t\ge 0\}$ of the sequence $\{Y_t^{(k)}: t\ge 0\}$ is the
solution of the martingale problem.
Let $Y_{s,k}^{-1}(u)$ and $b_k^{-1}(u)$ be defined as in Proposition \ref{p4.10}.
 For every $G\in C^2(\mbb{R})$ and $f\in C^1[0,1]$ we use (\ref{4.14})
 and It\^{o}'s formula to get
 \beqlb\label{4.18}
 G(\<Y_t^{(k)},f\>)
\ar=\ar
G(\<Y_0^{(k)},f\>)
+ \sigma_k\int_0^t ds \int_{\mbb{N}_+}\int_{[0,k]}\int_0^{kY_{s-}^{(k)}(1)}
\Big\{G\Big(\<Y_{s-}^{(k)},f\>\cr
 \ar\ar
 +\,k^{-1}(z-1)[f(1)-F^{(k)}(s,\theta,u)]\Big)- G(\<Y_{s-}^{(k)},f\>)\Big\}
 \bar{\pi}^{(k)}(dz,d\theta )du \cr
 \ar\ar
 + \sigma_k \int_0^t ds\int_0^{b_k(k)}d\theta\int_0^{kY_{s-}^{(k)}(1)}
\Big\{G\Big(\<Y_{s-}^{(k)},f\> -k^{-1}[f(1)-\tilde{F}^{(k)}(s,\theta,u)]\Big)\cr
 \ar\ar
-\, G(\<Y_{s-}^{(k)},f\>)\Big\} du
 +\sigma_k \int_0^t ds\int_{b_k(k)}^{b_k(0)}d\theta\int_0^{kY_{s-}^{(k)}(1)}
\Big\{G\Big(\<Y_{s-}^{(k)},f\>\cr
 \ar\ar
 +\, k^{-1}\tilde{F}^{(k)}(s,\theta,u)\Big)-\, G(\<Y_{s-}^{(k)},f\>)\Big\} du
+\mbox{local mart.}\cr
 \ar=\ar
G(\<Y_0^{(k)},f\>)
+ \sigma_k\int_0^t ds \int_{\mbb{N}_+}\int_{[0,k]}\int_0^{kY_{s-}^{(k)}(1)}
\Big\{G\Big(\<Y_{s-}^{(k)},f\>\cr
 \ar\ar
 +\,k^{-1}(z-1)f(Y_{s,k}^{-1}(u/k)\vee (\theta/k))\Big)
 -\, G(\<Y_{s-}^{(k)},f\>)\Big\} \bar{\pi}^{(k)}(dz,d\theta )du \cr
 \ar\ar
 + \sigma_k \int_0^t ds\int_0^{b_k(k)}d\theta\int_0^{kY_{s-}^{(k)}(1)}
\Big\{G\Big(\<Y_{s-}^{(k)},f\> -k^{-1}f(Y_{s,k}^{-1}(u/k))\Big)\cr
 \ar\ar
-\, G(\<Y_{s-}^{(k)},f\>)\Big\} du
 +\sigma_k \int_0^t ds\int_{b_k(k)}^{b_k(0)}d\theta\int_0^{kY_{s-}^{(k)}(1)}
\Big\{G\Big(\<Y_{s-}^{(k)},f\>\cr
 \ar\ar
 +\, k^{-1}[f(b_k^{-1}(\theta)/k)-f(Y_{s,k}^{-1}(u/k))]
 1_{\{Y_{s,k}^{-1}(u/k)\leq  b_k^{-1}(\theta)/k\}}\Big)\cr
 \ar\ar
 -\, G(\<Y_{s-}^{(k)},f\>)\Big\} du
+\, \mbox{local mart.}\cr
\ar=\ar
G(\<Y_0^{(k)},f\>)
+ k\sigma_k\int_0^t ds \int_{[0,1]}Y_{s-}^{(k)}(dx)\int_{\mbb{N}_+}\int_{[0,1]}
\Big\{G\Big(\<Y_{s-}^{(k)},f\>\cr
 \ar\ar
 +\,k^{-1}(z-1)f(x \vee \theta)\Big)
 -\, G(\<Y_{s-}^{(k)},f\>)\Big\} \bar{\pi}^{(k)}(dz,kd\theta ) \cr
 \ar\ar
 + k\sigma_k b_k(k) \int_0^t ds\int_{[0,1]}
\Big\{G\Big(\<Y_{s-}^{(k)},f\> -k^{-1}f(x)\Big)
-G(\<Y_{s-}^{(k)},f\>)\Big\} Y_{s-}^{(k)}(dx)\cr
\ar\ar
 + k\sigma_k \int_0^t ds \int_{[0,1]}Y_{s-}^{(k)}(dx)\int_1^x
\Big\{G\Big(\<Y_{s-}^{(k)},f\> +\, k^{-1}[f(\theta)-f(x)]\Big)\cr
\ar\ar
-\, G(\<Y_{s-}^{(k)},f\>)\Big\} b_k(kd\theta)
+\, \mbox{local mart.}\cr
\ar=\ar
G(\<Y_0^{(k)},f\>)
+ k\sigma_k\int_0^t ds \int_{[0,1]}Y_{s-}^{(k)}(dx)\int_{\mbb{N}}\int_{[0,1]}
\Big\{G\Big(\<Y_{s-}^{(k)},f\>\cr
 \ar\ar
 +\,k^{-1}(z-1)f(x \vee \theta)\Big)
 -\, G(\<Y_{s-}^{(k)},f\>)\Big\} \bar{\pi}^{(k)}(dz,kd\theta ) \cr
\ar\ar
 + k\sigma_k \int_0^t ds \int_{[0,1]}Y_{s-}^{(k)}(dx)\int_{\{0\}}\int_x^1
\epsilon_k(s,x,\theta)
\bar{\pi}^{(k)}(dz,kd\theta )\cr
\ar\ar
+\, \mbox{local mart.},
 \eeqlb
where
 \beqnn
\epsilon_k(s,\theta,x)
\ar=\ar
\Big\{G\Big(\<Y_{s-}^{(k)},f\> -k^{-1}f(x)\Big)-G\Big(\<Y_{s-}^{(k)},f\> -k^{-1}f(\theta)\Big)\Big\}\cr
\ar\ar
-\Big\{G\Big(\<Y_{s-}^{(k)},f\>+k^{-1}[f(\theta)-f(x)]\Big)-G(\<Y_{s-}^{(k)},f\>)\Big\}.
 \eeqnn
It is elementary to see that
 \beqnn
k\sigma_k\int_{\{0\}}\int_x^1 \epsilon_k(s,x,\theta)\bar{\pi}^{(k)}(dz,kd\theta )
 \eeqnn
tends to zero uniformly as $k\to\infty$.
Let $G(x)=e^{-x}$, by letting $k\to \infty$ in (\ref{4.18}) we get (\ref{4.17})
for $f\in C^1[0,1]$. A simple approximation shows
the martingale problem (\ref{4.17}) actually holds for any $f\in C[0,1]$.
By the proof of Theorem 7.13 in Li (2011) we get the result.\qed

 Let
$\{0\leq a_1<a_2<\cdots<a_n=1\}$ be an ordered set of constants. Denote by
$\{Y_{t,a_i}:t\geq 0\}$ and $\{Y_{t,a_i}^{(k)}:t\geq 0\}$ the restriction of $\{Y_t:t\geq 0\}$ and $\{Y_t^{(k)}:t\geq 0\}$ to $[0,a_i]$, respectively.  Let
$Y_{t}(a_i):=Y_t[0,a_i]$ and $Y_{t}^{(k)}(a_i):=Y_t^{(k)}[0,a_i] $ for every $t\geq 0$, $i=1,2,\cdots,n$. By arguments similar to those in He and Ma (2012) we get following results.

\btheorem\label{t3.3} Suppose that Condition (5.A) is satisfied and $\sup_{k\geq 1}\sigma_k b_k(0)<\infty$.
If $Y_{0}^{(k)}$
converges weakly to $Y_{0}\in M[0,1]$, then
$\{(Y_{t,a_1}^{(k)}, \cdots, Y_{t,a_n}^{(k)}):t\geq 0\}$ converges to
$\{(Y_{t,a_1}, \cdots, Y_{t,a_n}):t\geq 0\}$ in distribution on
$D([0,\infty),M[0,a_1]\times \cdots\times M[0,a_n])$.\etheorem

\bcorollary\label{c3.4} Suppose that Condition (5.A) is satisfied and
$\sup_{k\geq 1}\sigma_k b_k(0)<\infty$. If  $(Y_{0}^{(k)}(a_1), \cdots, Y_{0}^{(k)}(a_n))$
converges to $(Y_{0}(a_1), \cdots, Y_{0}(a_n))$, then $\{(Y_{t}^{(k)}(a_1), \cdots,
Y_{t}^{(k)}(a_n)):t\geq 0\}$ converges to $\{(Y_{t}(a_1), \cdots,
Y_{t}(a_n)):t\geq 0\}$ in distribution on $D([0,\infty), \mbb{R}_+^{n})$.
\ecorollary

\textbf{Acknowledgement.}  We would like to give our sincere thanks to
 Professor Zenghu Li for his encouragement and helpful discussions.


\bigskip

\bigskip

\noindent{\small Hui He and Rugang Ma: Laboratory of Mathematics and
Complex Systems, School of Mathematical Sciences, Beijing Normal
University, Beijing 100875, People's Republic of China. \\
\textit{E-mail:}{\tt hehui@bnu.edu.cn, marugang@mail.bnu.edu.cn} }

\end{document}